\newif\ifims 
\documentclass[11pt]{article}
\RequirePackage{amsthm,amsmath,amssymb,array,fancyhdr}
\RequirePackage[numbers]{natbib}
\usepackage[usenames]{color}
\RequirePackage[colorlinks,citecolor=blue,urlcolor=blue]{hyperref}
\usepackage[dvips]{graphicx}
\graphicspath{{./}{eps/}}
\newcommand{\bbC}{{\mathbb{C}}}
\newcommand{\bbN}{{\mathbb{N}}}
\newcommand{\bbR}{{\mathbb{R}}}
\newcommand{\bbZ}{{\mathbb{Z}}}
\newcommand{\cB}{{\mathcal{B}} }
\newcommand{\cE}{{\mathcal{E}}}
\newcommand{\cF}{{\mathcal{F}}}
\newcommand{\cN}{{\mathcal{N}}}
\newcommand{\fA}{{\mathfrak{A}}}
\renewcommand{\P}{{\textsf{P}}} \newcommand{\E}{{\textsf{E}}}
\newcommand{\Corr}{{\textsf{Corr}}}
\newcommand{\ind}{\mathrel{\mathop{\sim}\limits^{\mathrm{ind}}}}
\newcommand{\iid}{\mathrel{\mathop{\sim}\limits^{\mathrm{iid}}}}
\newcommand{\Bi}{\textsf{Bi}}
\newcommand{\NB}{\textsf{NB}}\newcommand{\No}{\textsf{No}}
\newcommand{\Po}{\textsf{Po}}
\newcommand{\AR}[1]{\textsf{AR}(#1)}
\newtheorem{thm}{Theorem}

\newcommand{\eps}{\epsilon}\newcommand{\hide}[1]{}

\newcommand{\set}[1]{\left\{#1\right\}}
\newcommand{\bet}[1]{\left[#1\right]}
\newcommand{\Bet}[1]{\Big[#1\Big]}
\newcommand{\cet}[1]{\left(#1\right)}

\newcommand{\ie}{\textit{i.e.{}}}
\newcommand{\pg}{\textit{p.}\thinspace}

\newcommand{\Sec}[1]{Sec.\thinspace(\ref{#1})}
\newcommand{\Secs}[2]{Secs.\thinspace\ref{#1} and \ref{#2}}
\newcommand{\Thm}[1]{Theorem\thinspace\ref{#1}}

\newcommand{\Eqn}[1]{Eqn.\thinspace(\ref{#1})}
\newcommand{\Eqns}[2]{Eqns.\thinspace(\ref{#1},\thinspace\ref{#2})}
\newcommand{\Eqss}[2]{Eqns.\thinspace(\ref{#1}--\ref{#2})}
\newcommand{\Fig}[1]{Fig.\thinspace(\ref{#1})}

\newcommand{\NNN}[3]{{N_{#1#2#3}}}
\renewcommand{\lll}[3]{{\lambda_{\,#1#2#3}}}

\newcommand{\rr}[1]{{r_{#1}}}
\renewcommand{\th}[1]{{\theta_{#1}}}

\newcommand{\omp}{(1{-}p)}

\newcommand{\omrh}{(1{-}\rho)}

\newcommand{\jgz}{_{j\ge1}}            
\newcommand{\nljgz}{\nolimits_{j\ge1}} 
\newcommand{\iz}{_{i\ge0}}             
\newcommand{\nliz}{\nolimits_{i\ge0}}  
\newcommand{\kz}{_{k\ge0}}             
\newcommand{\nlkz}{\nolimits_{k\ge0}}  
\newcommand{\ikz}{_{i,k\ge0}}          
\newcommand{\nlikz}{\nolimits_{i,k\ge0}} 
\newcommand{\zz}[1]{\ifcase #1{*}\or{s}\or{z}\or{t}\else{*}\fi}
\newcommand{\vpz}[1]{{\varphi(\zz1,\zz3\mid #1)}}

\newcommand{\T}{T} 
\newcommand{\sT}{{\scriptscriptstyle T}} 
\newcommand{\dT}{{|T|}} 
\newcommand{\mT}{{\mu_\sT}}
\newcommand{\ott}{_{\{1,2,3\}}}
\makeatletter 
\multiply\count2 60 \advance\count2 \time
\edef\now{\two@digits{\the\count1}:\two@digits{\the\count2}}
\ifx\newcolumntype\undefined
\else
  \newcolumntype{C}{>{$}c<{$}}
  \newcolumntype{L}{>{$}l<{$}}
  \newcolumntype{R}{>{$}r<{$}}
\fi
\def\rev$Revi#1: #2 ${\gdef\vsn{#2}}

\newcommand{\spt}{\mathop{\mathrm{spt}}}
\newcommand{\BB}{{\mathsf{BB}}}
\newcommand{\bbU}{{\mathbb{U}}}
\rev$Revision: 1.1 $
\newif{\ifBib} \Bibfalse
\begin{document}
\title
{Markov Infinitely-Divisible Stationary\\ Time-Reversible Integer-Valued
      Processes}
\author{Robert L. Wolpert$^{a*}$~\& Lawrence D. Brown$^b$\\
  { $^a$Department of Statistical Science}\\
  { Duke University, Durham NC 27708-0251, USA}\\
  { $^b$Deceased, 2017-02-21}\\ {}}
\date{\today} 
\maketitle
\renewcommand{\abstractname}{Summary}
\begin{abstract}
  We prove a complete class theorem that characterizes \emph{all}
  stationary time reversible Markov processes whose finite dimensional
  marginal distributions (of all orders) are infinitely divisible.  Aside
  from two degenerate cases (iid and constant), in both discrete and
  continuous time every such process with full support is a branching
  process with Poisson or Negative Binomial marginal univariate
  distributions and a specific bivariate distribution at pairs of times.
  As a corollary, we prove that every nondegenerate stationary integer
  valued processes constructed by the Markov thinning process fails to have
  infinitely divisible multivariate marginal distributions, except for the
  Poisson.  These results offer guidance to anyone modeling integer-valued
  Markov data exhibiting autocorrelation.
\end{abstract}

\section{Introduction}\label{s:intro}

Many applications feature autocorrelated count data $X_t$ at discrete times
$t$.  A number of authors have constructed and studied stationary stochastic
processes $X_t$ whose one-dimensional marginal distributions come from an
arbitrary infinitely-divisible distribution family $\{\mu^{\theta}\}$, such
as the Poisson $\Po(\theta)$ or negative binomial $\NB(\theta,p)$, and that
are ``$\AR1$-like'' in the sense that their autocorrelation function is
$\Corr[X_s,X_t] =\rho^{|s-t|}$ for some $\rho\in (0,1)$ \citep {Lewi:1983,
  Lewi:McKe:Hugu:1989, McKe:1988, AlOs:Alza:1987, Joe:1996}.  The most common
approach is to build a time-reversible Markov process using \emph{thinning},
in which the process at any two consecutive times may be written in the form
\[ X_{t-1}=\xi_t+\eta_t \qquad X_t=\xi_t+\zeta_t \] with $\xi_t$, $\eta_t$,
and $\zeta_t$ all independent and from the same infinitely-divisible family
(see \Sec{ss:thin} below for details).  A second construction of a stationary
time-reversible process with the same one-dimensional marginal distributions
and autocorrelation function, with the feature that its finite-dimensional
marginal distributions of all orders are infinitely-divisible, is to set
$X_t:=\cN(G_t)$ for a \emph {random measure} $\cN$
\iftrue 
on some measure space $(E,\cE,m)$ that assigns independent
infinitely-divisible random variables $\cN(A_i)\sim\mu^{\theta_i}$ to
disjoint sets $A_i\in\cE$ of measure $\theta_i=m(A_i)$, and a family of sets
$\{G_t\}\subset\cE$ whose intersections have measure $m\big(G_s\cap
G_t\big)=\theta \rho^{|s-t|}$ \else on some space that assigns independent
infinitely-divisible random variables $\cN(A_i)\sim\mu^{|A_i|}$ to disjoint
sets $A_i$ of measure $|A_i|$, and a family of sets $\{G_t\}$ whose
intersections have measure $|G_s\cap G_t| =\theta \rho^{|s-t|}$ \fi (see
\Sec{ss:meas}).

For the normal distribution $X_t\sim\No(\mu,\sigma^2)$, these two
constructions both yield the usual Gaussian $\AR1$ process.  The two
constructions also yield identical processes for the Poisson $X_t\sim
\Po(\theta)$ distribution, but they differ for all other nonnegative
integer-valued infinitely-divisible distributions.  For each nonnegative
integer-valued infinitely-divisible marginal distribution except the Poisson,
the process constructed by thinning does not have infinitely-divisible
marginal distributions of all orders (\Thm{t:thin}, \Sec{ss:thm}), and the
process constructed using random measures does not have the Markov property
(\Thm{t:meas}, \Sec{ss:thm}).  Thus none of these is completely satisfactory
for modeling autocorrelated count data with heavier tails than the Poisson
distribution.

In the present manuscript we construct and characterize every process that is
Markov, infinitely-divisible, stationary, and time-reversible with
non-negative integer values.  The formal characterization is contained in the
statement of \Thm{t:struc} in \Sec {ss:thm}, which follows necessary
definitions and the investigation of special cases needed to establish the
general result.

\subsection{Thinning Process}\label{ss:thin}

Any univariate infinitely-divisible (ID) distribution $\mu(dx)$ on $\bbR^1$
is $\mu^1$ the for a convolution semigroup $\{\mu^\theta:~\theta\ge0\}$ and,
for $0<\theta<\infty$ and $0<\rho<1$, determines uniquely a ``thinning
distribution'' $\mu^\theta_\rho(dy\mid x)$ of $Y$ conditional on the sum
$X=Y+Z$ of independent $Y\sim\mu^{\rho\theta}$ and $Z\sim\mu^{\omrh \theta}$.
This thinning distribution determines a unique stationary time-reversible
Markov process with one-step transition probability distribution given by the
convolution
\[ \P[X_{t+1}\in A\mid \cF_t ] 
   = \int_{\xi+\zeta\in A} 
     \mu^{\omrh\theta}(d\zeta)~\mu^\theta_\rho (d\xi\mid X_t) 
\]
for Borel sets $A\subset\bbR$, where $\cF_t=\sigma\{X_s:~s\le t\}$ is the
minimal filtration.  By induction the auto-correlation is $\Corr (X_s,X_t)
=\rho^{|s-t|}$ for square-integrable $X_t$.  The process can be constructed
beginning at any $t_0\in\bbZ$ by setting
\begin{subequations}\label{e:thin}
\begin{align}
X_{t_0}&\sim \mu^\theta(dx)\label{e:thin1}\\
\xi_t &\sim  \mu^\theta_\rho (d\xi\mid x) \text{ with }
       x=\begin{cases} X_{t-1}& \text{if }t>t_0\\
                       X_{t+1}&\text{if } t<t_0\end{cases}\label{e:thin2}\\
X_t &:= \xi_t+\zeta_t\qquad\text{for }
        \zeta_t\sim \mu^{\theta\omrh}(d\zeta).\label{e:thin3}
\end{align}
\end{subequations}
Time-reversibility and hence the lack of dependence of this definition on the
choice of $t_0$ follows from the argument presented in the proof of
\Thm{t:thin} in \Sec{ss:thm} below.

\subsubsection{Thinning Example 1: Poisson}\label{sss:pth}
For Poisson-distributed $X_t\sim\mu^\theta=\Po(\theta)$ with mean $\theta>0$,
for example, the thinning recursion step for $0<\rho<1$ and $t>t_0$ can be
written
\begin{align*}
  X_{t} &= \xi_t + \zeta_t \text{\quad for independent:}\\
  \xi_t &\sim \Bi\big(X_{t-1}, \rho\big), \qquad
\zeta_t  \sim \Po\big(\theta\omrh\big)
\end{align*}
and hence the joint generating function at two consecutive times is
\begin{alignat*}3
\phi(s,z) &= \E\Big[ s^{X_{t-1}} z^{X_t}\Big] &
          &= \exp\Big\{(s+z-2)\theta\omrh + (s\,z-1)\theta\rho\Big\}
. 
\end{alignat*}
This was called the ``Poisson $\mathrm{AR}(1)$ Process'' by 
McKenzie \citep{McKe:1985} and has been studied by many other
authors since its introduction.

\subsubsection{Thinning Example 2: Negative Binomial}
In the thinning process applied to the Negative Binomial $X_t\sim \mu^\theta
=\NB(\theta,p)$ distribution with mean $\theta\omp/p$, recursion for $t>t_0$
takes the form
\begin{align}
  X_{t} &= \xi_t + \zeta_t \text{\quad for independent:}\notag\\
  \xi_t &\sim \BB\big(X_{t-1};~ \theta\rho,~\theta\omrh\big), \qquad
  \zeta_t \sim \NB\big(\theta\omrh,~p\big)\notag
\intertext{for beta-binomial distributed $\xi_t \sim \BB
 (n;\alpha,\beta)$ \citep[see] [\S2.2] {John:Kemp:Kotz:2005} with
 $n=X_{t-1}$, $\alpha=\theta \rho$, and $\beta= \theta \omrh$, and
 negative binomial $\zeta_t \sim \NB\big(\theta \omrh,p\big)$.  Thus the
 joint generating function is}
  \phi(s,z) &= \E\Big[ s^{X_{t-1}} z^{X_t}\Big]\notag\\
  &= p^{\theta(2-\rho)} (1-q\,s)^{-\theta\omrh}\, (1-q\,z)^{-\theta\omrh}\,
  (1-q\,s\,z)^{-\theta \rho}.\label{e:phi-sz-nbt}
\end{align}
From this one can compute the conditional generating function
\[
\phi(z\mid x)=\E\bet{z^{X_{t}}\mid X_{t-1}=x}
              =  \cet{\frac{p}{1-qz}}^{\theta\omrh}~
                 {}_2F_1(\theta\rho, -x; \theta; 1-z)
\]
where $_2F_1(a,b;c;z)$ denotes Gauss' hypergeometric function \citep [\S15]
{Abra:Steg:1964} and, from this (for comparison below),
\begin{align}
\P[ X_{t-1}=0, X_{t+1}=0\mid X_t=2] &=
   [p^{\theta\omrh}~{}_2F_1(\theta\rho, -x; \theta; 1)]^2\notag\\
   &= \bet{p^{\theta\omrh} \omrh}^2
      \bet{\frac{1+\theta\omrh}{1+\theta}}^2.\label{e:nbt020}
\end{align}
This process, as we will see below in \Thm{t:thin}, is Markov, stationary,
and time-reversible, with infinitely-divisible one-dimensional marginal
distributions $X_t\sim\NB(\theta,p)$, but the joint marginal distributions at
three or more consecutive times are not ID.  It appears to have been
introduced by Joe 
\citep[\pg665] {Joe:1996}.

\subsection{Random Measure Process}\label{ss:meas}
Another approach to the construction of processes with specified univariate
marginal stationary distribution $\mu^\theta(dx)$ is to set $X_t := \cN(G_t)$
for a \emph {random measure} $\cN$ and a class of sets $\set{G_t}$,
as in \citep [\S3.3, 4.4] {Wolp:Taqq:2005}.  We begin with a countably
additive random measure $\cN(dx\,dy)$ that assigns independent random
variables $\cN(A_i)\sim\mu^{|A_i|}$ to disjoint Borel sets $A_i\in\cB
(\bbR^2)$ of finite area $|A_i|$ (this is possible by the Kolmogorov
consistency conditions), and a collection of sets
\[ G_t :=\set{(x,y):~x\in\bbR, ~0 \le y < \theta\lambda~e^{-2\lambda
    |t-x|}}\]
\begin{figure}[!htp]
\begin{center}
\includegraphics[height=70mm]{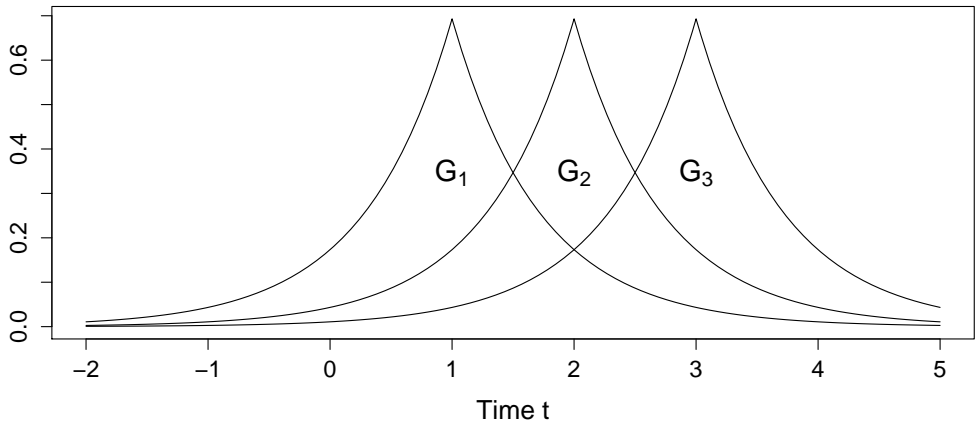}
\vspace*{-5mm}
\caption{\label{f:Gt}Random measure construction of process $X_t=\cN(G_t)$}
\end{center}
\vspace*{-5mm}
\end{figure}
(shown in \Fig{f:Gt}) whose intersections satisfy $|G_s\cap G_t| =
\theta e^{-\lambda|s-t|}$.  For $t\in\bbZ$, set
\iffalse
\begin{subequations}\label{e:meas}
\begin{align}
 X_t &:= \cN(G_t)\label{e:meas.X}
\intertext{for the set}
 G_t &:=\set{(x,y):~x\in\bbR, ~0 \le y < \theta\lambda~e^{-2\lambda
    |t-x|}}\label{e:meas.G}
\end{align}\end{subequations}
\else
\begin{equation}\label{e:meas}
 X_t := \cN(G_t).
\end{equation}
\fi 
For any $n$ times $t_1<t_2<\dots<t_n$ the sets $\set {G_{t_i}}$ partition
$\bbR^2$ into $n(n+1)/2$ sets of finite area (and one with infinite area,
$(\cup G_{t_i})^c$), so each $X_{t_i}$ can be written as the sum of some
subset of $n(n+1)/2$ independent random variables.  In particular, any $n=2$
variables $X_s$ and $X_t$ can be written as
\[ X_s=\cN(G_s\backslash G_t)+\cN(G_s \cap G_t),\qquad
   X_t=\cN(G_t\backslash G_s)+\cN(G_s \cap G_t)
\]
just as in the thinning approach, so both 1-dimensional and 2-dimensional
marginal distributions for the random measure process coincide with those
for the thinning process of \Sec{ss:thin}.

Evidently the process $X_t$ constructed from this random measure is
stationary, time-reversible and infinitely divisible in the strong sense that
all finite-dimensional marginal distributions are ID.  Although the 1- and
2-dimensional marginal distributions of this process coincide with those of
the thinning process, the $k$-dimensional marginals may differ for $k\ge3$,
so this process cannot be Markov.  We will see in \Thm{t:meas} below that the
only nonnegative integer-valued distribution for which it is Markov is the
Poisson.

\subsubsection{Random Measure Example 1: Poisson}\label{sss:prm}
The conditional distribution of $X_{t_n}=\cN(G_{t_n})$ given $\{X_{t_j}:
~j<n\}$ can be written as the sum of $n$ independent terms, $(n-1)$ of them
with binomial distributions (all with the same probability parameter
$p=\rho^{|t_n-t_{n-1}|}$, and with size parameters that sum to $X_{t_{n-1}}$)
and one with a Poisson distribution (with mean $\theta(1-\rho^
{|t_n-t_{n-1}|}$).  It follows by induction that the random-measure Poisson
process is identical in distribution to the thinning Poisson process of
\Sec{sss:pth}.

\subsubsection{Random Measure Example 2: Negative Binomial}
The random variables $X_1$, $X_2$, $X_3$ for the random measure process built
on the Negative Binomial distribution $X_t\sim\NB(\theta,p)$ with
autocorrelation $\rho\in(0,1)$ can be written as sums 
\[ X_1 = \zeta_1+\zeta_{12}+\zeta_{123} \qquad
   X_2 = \zeta_2+\zeta_{12}+\zeta_{23}+\zeta_{123} \qquad
   X_3 = \zeta_3+\zeta_{23}+\zeta_{123}
\]
of six independent negative binomial random variables
$\zeta_s\sim\NB(\theta_s, p)$ with shape parameters 
\[ \theta_1=\theta_3=\theta\omrh,\qquad
   \theta_2=\theta\omrh^2,\qquad
   \theta_{12}=\theta_{23}=\theta\rho\omrh,\qquad
   \theta_{123}=\theta\rho^2
\]
(each $\zeta_s=\cN\big(\cap_{t\in s} G_t\big)$ and $\theta_s=|\cap_{t\in s}
G_t|$ in \Fig{f:Gt}).  It follows that the conditional probability
\begin{align}
\P[ X_1=0,~X_3=0\mid X_2=2] 
  &= \P[\zeta_1=\zeta_{12}=\zeta_{123}=\zeta_{23}=\zeta_3=0 \mid
        \zeta_2+\zeta_{12}+\zeta_{23}+\zeta_{123}=2]\notag\\
  &=\frac{\P[\zeta_2=2, \text{ all other } \zeta_s=0]} {\P[X_2=2]}\notag\\
  &= \bet{p^{\theta\omrh}\omrh}^2{\frac{1+\theta\omrh^2}{1+\theta}}\label{e:nbrm020}
\end{align}
differs from that of the thinning negative binomial process in \Eqn{e:nbt020}
for all $\theta>0$ and $\rho>0$.  Thus this process is stationary,
time-reversible, and has infinitely-divisible marginal distributions of all
orders, but it cannot be Markov since its 2-dimensional marginal
distributions coincide with those of the Markov thinning process but its
3-dimensional marginal distributions do not.

In \Sec {s:2sol} of this paper we characterize every discrete-time process
that is Markov, Infinitely-divisible, Stationary, and Time-reversible with
non-negative Integer values (\emph{MISTI} for short).  In \Sec {s:simm} we
first present the necessary definitions and preliminary results; in \Sec
{s:cont} we extend the results to continuous time, with discussion in \Sec
{s:disc}.

\section{MISTI Processes}\label{s:simm}
A real-valued stochastic process $X_t$ indexed by $t\in\bbZ$ is
\emph{stationary} if each finite-dimensional marginal distribution
 \[ \mT(B):=\P\bet{ X_\T\in B} \]
\begin{subequations}\label{e:props}%
satisfies 
\begin{equation}
\mT(B)=\mu_{s+\sT}(B)\label{e:stat}
\end{equation}
for each set $\T\subset\bbZ$ of cardinality $|\T|<\infty$, Borel set
$B\in\cB(\bbR^\dT)$, and $s\in\bbZ$, where as usual ``$s+T$'' denotes
$\{(s+t):~t\in\T\}$.  A stationary process is \emph{time-reversible} if also
\begin{equation}
\mT(B)=\mu_{-\sT}(B)\label{e:trev}
\end{equation}
(where ``$-T$'' is $\{-t:~t\in T\}$) and \emph{Markov} if for every
$t\in\bbZ$ and finite $\T\subset\{s\in\bbZ:~s\ge t\}$,
\begin{equation}
        \P[ X_\T\in B\mid \cF_t] = \P[ X_\T\in B\mid X_t] \label{e:mark}
\end{equation}
for all $B\in\cB(\bbR^\dT)$, where $\cF_t:=\sigma\{X_s:~s\le t\}$.  The
process $X_t$ is Infinitely Divisible (ID) or, more specifically,
\emph{multivariate} infinitely divisible (MVID) if each $\mT$ is the $n$-fold
convolution of some other distribution $\mu_\sT^{(1/n)}$ for each $n\in\bbN$.
This is more restrictive than requiring only that the one-dimensional
marginal distributions be ID and, for integer-valued processes that satisfy
\begin{equation}
 \mT(\bbZ^\dT)=1, \label{e:integ}
\end{equation}
it is equivalent by the L\'evy-Khinchine formula \cite
[\pg74] {roge:will:2000a} to the condition that each $\mT$ have
characteristic function of the form
\begin{equation}
\int_{\bbR^\dT} e^{i\omega'x}\,\mT(dx)
         = \exp\set{\int_{\bbZ^\dT} \big(e^{i\omega'u}-1\big)\,
               \nu_\sT(du)},\qquad \omega\in\bbR^\dT \label{e:mvid}
\end{equation}
\end{subequations}
for some finite measure $\nu_\sT$ on $\cB(\bbZ^\dT)$.  Call a process $X_t$
or its distributions $\mT(du)$ \emph{MISTI} if it is 
Markov, 
nonnegative Integer-valued, 
Stationary,
Time-reversible, 
and 
Infinitely divisible,  
\ie, satisfies \Eqss {e:stat} {e:mvid}.  We now turn to the problem of
characterizing all MISTI distributions.

\subsection{Three-dimensional Marginals}\label{ss:3dim}

By stationarity and the Markov property all MISTI finite-dimensional
distributions $\mT(du)$ are determined completely by the marginal
distribution for $X_t$ at two consecutive times; to exploit the MVID property
we will study the three-dimensional marginal distribution for $X_t$ at any
set $T$ of $\dT=3$ consecutive times--- say, $T=\{1,2,3\}$.  By \Eqn{e:mvid}
we can represent $X\ott$ in the form
\[
   X_1 = \sum i \NNN i++\qquad
   X_2 = \sum j \NNN +j+\qquad
   X_3 = \sum k \NNN ++k
\]
for independent Poisson-distributed random variables
\[ \NNN ijk \ind \Po(\lll ijk) \] with means $\lll ijk := \nu(\{(i,j,k)\})$;
here and hereafter, a subscript ``$+$'' indicates summation over the entire
range of that index--- $\bbN_0=\set{0,2\dots}$ for $\set{\NNN ijk}$ and
$\set{\lll ijk}$, $\bbN=\set{1,2,\dots}$ for $\set{\th j}$.  The sums $\th
j:=\lll+j+$ for $j\ge1$ characterize the univariate marginal distribution of
each $X_t$--- for example, through the probability generating function (pgf)
\[ \varphi(\zz2) := \E[ \zz2^{X_t} ] 
      = \exp\left[\sum\nljgz \big(\zz2^j-1\big)\th j\right].
\]
To avoid trivial technicalities we will assume that $0 < \P[X_t=1]=\varphi'(0)
  =\th1e^{-\th+}$, \ie, $\th1>0$.  Now set $\rr i:=\lll i1+/\th1$, and for
  later use define functions: 
\begin{equation}\label{e:pP}
 \psi_j(\zz1,\zz3) := \sum\ikz \zz1^i\zz3^k\lll ijk\qquad\qquad
     p(\zz1) := \psi_1(s,1)/\th1=\sum\iz \zz1^i \,\rr i\qquad\qquad
     P(\zz2) := \sum\jgz \zz2^j \,\th j.
\end{equation}
Since $\rr i$ and $\th j$ are nonnegative and summable (by \Eqns
{e:integ}{e:mvid}), $p(\zz1)$ and $P(\zz2)$ are analytic on the open unit
ball $\bbU \subset \bbC$ and continuous on its closure.  Similarly, since $\lll
ijk$ is summable, each $\psi_j(\zz1,\zz3)$ is analytic on $\bbU^2$ and
continuous on its closure.  Note $\psi_j(1,1)=\th j$, $p(0)=\rr0$ and
$p(1)=1$, while $P(0)=0$ and $P(1)=\th+$; also $\varphi(\zz2) =
\exp\set{P(\zz2)-\th+}$.  Each $\psi_j(s,t)=\psi_j(t,s)$ is symmetric by
\Eqn{e:trev}, as are the conditional probability generating functions:
\begin{equation}\notag 
    \vpz z := \E\big[\zz1^{X_1}\zz3^{X_3}\mid X_2=z\big].
\end{equation}

\subsubsection{Conditioning on $X_2=0$}\label{sss:x2=0}
By the Markov property \Eqn{e:mark}, $X_1$ and $X_3$ must be conditionally
independent given $X_2$, so the conditional probability generating function
must factor:
\begin{align}
    \vpz0 &:= \E\big[\zz1^{X_1}\zz3^{X_3}\mid X_2=0\big] 
            = \E\big[\zz1^{\sum\iz i\NNN i0+}~
                    \zz3^{\sum\kz k\NNN+0k}\big]\notag\\
          &= \exp\Big\{\sum\ikz
                 \big(\zz1^i\zz3^k-1\big)\lll i0k\Big\}\notag\\
          &\equiv \varphi(\zz1,1\mid0)~ \varphi(1,\zz3\mid0).\label{e:fac0}
\intertext{Taking logarithms,}
             {\sum\big(\zz1^i\zz3^k-1\big)\lll i0k} 
          &\equiv {\sum\big(\zz1^i-1\big)\lll i0k}
            +{\sum\big(\zz3^k-1\big)\lll i0k}\notag
\intertext{or, for all $\zz1$ and $\zz3$ in the unit ball in $\bbC$,}
         0&\equiv\sum(\zz1^i-1)(\zz3^k-1\big)\lll i0k.\label{e:x2=0}
\intertext{Thus $\lll i0k=0$ whenever both $i>0$ and $k>0$ and, by symmetry,}
    \varphi(1,z\mid0)=\varphi(z,1\mid0)
    &=\exp\set{\sum\nliz (z^i-1)\lll i00}.\notag
\end{align}

\subsubsection{Conditioning on $X_2=1$}\label{sss:x2=1}
Similarly
\begin{align}
    \vpz1 := \E\big[\zz1^{X_1}\zz3^{X_3}\mid X_2=1\big]
           &= \E\big[\zz1^{\sum\iz i(\NNN i0+ + \NNN i1+)}\quad
              \zz3^{\sum\kz k(\NNN +0k + \NNN+1k)} \mid
              \NNN+1+=1\big]\notag\\ 
           &= \vpz0 \E\big[\zz1^{\sum\iz i\NNN i1+}\quad
                          \zz3^{\sum\kz k\NNN +1k} \mid
                          \NNN+1+=1\big]\notag\\ 
          &= \vpz0 \set{\sum\ikz \zz1^i\zz3^k\,\big[\lll i1k/\lll
          +1+\big]}\notag 
\intertext{since $\set{\NNN i1k}$ is conditionally multinomial given $\NNN
          +1+$ and independent of $\set{\NNN i0k}$.  By the Markov property
          this too must factor, as
           $\varphi(\zz1,\zz3\mid1)=\varphi(\zz1,1\mid1)
          \,\varphi(1,\zz3\mid1)$, so by \Eqn{e:fac0}}
          \th1 \set{\sum\nlikz \zz1^i\zz3^k\lll i1k}
           &= \set{\sum\nliz \zz1^i \lll i1+}
              \set{\sum\nlkz \zz3^k \lll +1k}\notag
\intertext{or, since $\lll i1k=\lll k1i$ by \Eqns{e:trev}{e:pP},}
\psi_1(\zz1,\zz3) &:= \sum\nlikz \zz1^i\zz3^k \lll i1k
                  = \th1 p(\zz1)\,p(\zz3),\notag\\
           \vpz1 &= \vpz0\, p(\zz1)\, p(\zz3).\notag
\end{align}

\subsubsection{Conditioning on $X_2=2$}\label{sss:x2=2}
The event $\set{X_2=2}$ for $X_2:=\sum\jgz j\NNN+j+$ can happen in two ways:
either $\NNN+1+ =2$ and each $\NNN+j+=0$ for $j\ge2$, or $\NNN+2+=1$ and
$\NNN+j+=0$ for $j=1$ and $j\ge3$, with $\NNN+0+$ unrestricted in each case.
These two events have probabilities $(\th1^2/2)e^{-\th+}$ and
$(\th2)e^{-\th+}$, respectively, so the joint generating function for
$\{X_1,X_3\}$ given $X_2=2$ is
\begin{align}
    \vpz2 &:= \E\big[\zz1^{X_1}\zz3^{X_3}\mid X_2=2\big]\notag\\
          &= \E\big[\zz1^{\sum\iz i(\NNN i0+ + \NNN i1+ + \NNN i2+)}\quad
                    \zz3^{\sum\kz k(\NNN +0k + \NNN +1k + \NNN +2k)} 
                    \mid \NNN+1++2\NNN+2+=2\big]\notag\\
          &= \vpz0\set{
                   \frac{\th1^2/2}{\th1^2/2+\th2}
                   \bet{\sum\ikz \zz1^i\zz3^k\lll i1k/\lll+1+}^2  +
                   \frac{\th2}{\th1^2/2+\th2}
                   \bet{\sum\ikz \zz1^i\zz3^k\lll i2k/\lll+2+} }\notag\\
          &= \frac{\vpz0}{\th1^2/2+\th2}\set{
                   \frac{\th1^2}{2}
                   \bet{\sum\ikz \zz1^i\zz3^k\lll i1k/\th1}^2  +
                   \th2 \bet{\sum\ikz \zz1^i\zz3^k\lll i2k/\th2} }\notag\\
          &= \frac{\vpz0}{\th1^2/2+\th2}~
                   \set{\frac{\th1^2}2 p(\zz1)^2 p(\zz3)^2 +
                   \psi_2(s,t)}.\label{e:x2=2}
\end{align}
In view of \Eqn{e:fac0}, this will factor in the form $\varphi(\zz1,\zz3\mid2)=
  \varphi(\zz1,1\mid2)\, \varphi(1,\zz3\mid2)$ as required by Markov property
  \Eqn{e:mark} if and only if for all $s,t$ in the unit ball:
\begin{align}
     \bet{\frac{\th1^2}{2}+\th2}
     \bet{\frac{\th1^2}2 p(\zz1)^2 p(\zz3)^2 + \psi_2(s,t)}
 &=  \bet{\frac{\th1^2}2 p(\zz1)^2 + \psi_2(s,1)}
     \bet{\frac{\th1^2}2 p(\zz3)^2 + \psi_2(1,t)} \notag 
\end{align}
or
\begin{multline}
 \frac{\th1^2}2 \Bet{\th2 p(s)^2 p(t)^2 - p(s)^2 \psi_2(1,t)
             - \psi_2(s,1) p(t)^2 + \psi_2(s,t)} \\
    = \Bet{ \psi_2(s,1)\,\psi_2(1,t)-\th2\psi_2(s,t) }. \notag 
\end{multline}
To satisfy the ID requirement of \Eqn {e:mvid}, this must hold with each $\th
j$ replaced by $\th j/n$ for each integer $n\in\bbN$.  Since the left and
right sides are homogeneous in $\theta$ of degrees $3$ and $2$ respectively,
this will only happen if each square-bracketed term vanishes identically,
\ie, if
\begin{align}
 \th2 \psi_2(s,t) &\equiv  \psi_2(s,1)\psi_2(1,t)\notag
\intertext{and}
0 &= 
  \th2 \Bet{\th2 p(s)^2 p(t)^2 - p(s)^2 \psi_2(1,t)
             - \psi_2(s,1) p(t)^2} + \psi_2(s,1) \psi_2(1,t)\notag\\
  &= \bet{\th2 p(s)^2-\psi_2(s,1)}  \bet{\th2 p(t)^2-\psi_2(1,t)},\notag
\intertext{so}
 \psi_2(s,t) & :=  \sum\ikz \zz1^i\zz3^k\lll i2k
                = \th2 p(s)^2\,p(t)^2,\notag\\ 
        \vpz2   &= \vpz0\, p(\zz1)^2 p(\zz3)^2.\notag
\end{align}

\subsubsection{Conditioning on $X_2=j$}\label{sss:x2=j}
The same argument applied recursively, using the Markov property for each
$j\ge1$ in succession, leads to:
\begin{multline} \notag 
     \Bet{\frac{\th1^j}{j!}+\dots+\th1\th{j-1}}
     \Bet{\th j p(s)^j p(t)^j - p(s)^j\psi_j(1,t)
         -\psi_j(s,1)p(t)^j + \psi_j(s,t)}\\
 = \Bet{\psi_j(s,1)\psi_j(1,t)- \th j \psi_j(s,t)}
\end{multline}
so
\begin{equation}
 \psi_j(s,t) := \sum\ikz \zz1^i\zz3^k\lll ijk
              = \th j\, p(\zz1)^j p(\zz3)^j,\qquad j\ge1\label{e:key}
\end{equation}
and consequently
\begin{align}
    \vpz j &= \E\big[\zz1^{X_1}\zz3^{X_3}\mid X_2=j\big] 
            = \bet{\varphi(\zz1,1\mid0)\, 
              p(\zz1)^j} ~
              \bet{\varphi(1,\zz3\mid0)\,
              p(\zz3)^j}.\notag
\end{align}
Conditionally on $\set{X_2=j}$, $X_1$ and $X_3$ are distributed
independently, each as the sum of $j$ independent random variables with
generating function $p(s)$, plus one with generating function
$\varphi(s,1\mid0)$--- so $X_t$ is a branching process \citep {Harr:1963}
whose unconditional three-dimensional marginal distributions have generating
function:
\begin{align}
\varphi(\zz1,\zz2,\zz3) &:= \E\big[\zz1^{X_1}\zz2^{X_2}\zz3^{X_3}\big]\notag\\
       &= \vpz0 \sum_{j\ge0} \zz2^j p(\zz1)^j p(\zz3)^j \P[X_2=j]\notag\\
       &= \vpz0 \E\left[ z p(\zz1) p(\zz3)\right]^{X_2}\notag\\
       &= \vpz0 \varphi\big( z p(\zz1) p(\zz3)\big)\notag\\
       &= \vpz0 \exp\big[P\big( z p(\zz1) p(\zz3)\big)-\th+\big].\label{e:gen3}
\end{align}
See \Secs {ss:pops} {s:disc} for further development of this branching
process representation.

\subsection{Stationarity}\label{ss:sta}

Without loss of generality we may take $\lll 000 = 0$.  By \Eqn{e:key} with
$\zz1=0$ and $\zz3=1$ we have $\lll 0j+=\th j\rr0^j$; by \Eqn{e:x2=0} we have
$\lll i00 = \lll i0+$.  By time-reversibility we conclude that $\lll i00 = 0$
for $i=0$ and, for $i\ge1$,
\begin{equation}
  \lll i00 = \th i \rr0^i \label{e:i00}.
\end{equation}
Now we can evaluate
\[
  \varphi(\zz1,\zz3\mid0) = \exp\set{P(\zz1\,\rr0)+P(\zz3\,\rr0)-2P(\rr0)}
\]
and, from this and \Eqn{e:gen3}, evaluate the joint generating function for
$X_{\{1,2,3\}}$ as:
\begin{align}
 \varphi(\zz1,\zz2,\zz3) 
     &=\exp\set{P\big( z\, p(\zz1) p(\zz3)\big)-\th+
       +P(\zz1\,\rr0)+P(\zz3\,\rr0)-2P(\rr0)},\qquad j\ge1 \label{e:phi}
\intertext{and so that for $X_{\{1,2\}}$ as:} 
  \varphi(\zz1,\zz2,1) 
         &=\exp\set{P\big( z\, p(\zz1)\big)-\th+ 
                  + P(\zz1\,\rr0)-P(\rr0)}. \label{e:x1x2}
\end{align}
Now consider \Eqn{e:key} with $\zz3=1$,
\begin{align}
   \sum\iz \zz1^i\lll ij+ &= \th j\, p(\zz1)^j.\label{e:sz}
\intertext{It follows first for $j=1$ and then for $i=1$ that}
  \lll i1+ &= \th1 \rr i& i&\ge1\notag\\
  \lll 1j+ &= \th j [ j \rr0^{j-1} \rr1 ]&j&\ge1\notag
\intertext{so again by time reversibility with $i=j$, since $\th1>0$, we have}
\rr j &= \th j [j\, \rr0^{j-1} \rr1]/ \th1\qquad j\ge1. \label{e:rj}
\end{align}
Thus $\rr0$, $\rr1$, and $\{\th j\}$ determine all the $\{\rr j\}$ and so all
the $\{\lll ijk\}$ by \Eqns {e:key} {e:i00} and hence the joint distribution
of $\{X_t\}$.

Now consider \Eqn {e:sz} first for $j=2$ and then $i=2$:
\begin{align}
  \sum\iz \zz1^i \lll ij+ &= \th j\,\left[\sum\nliz \zz1^i \rr
  i\right]^j\notag \\
  \lll i2+ &= \th2 \sum_{k=0}^i \rr k \rr {i-k} & i&\ge2\notag\\
  \lll 2j+ &= \th j \left[
        j            \rr0^{j-1} \rr2   + 
        \binom{j}{2} \rr0^{j-2} \rr1^2 \right]  & j&\ge2\notag
\end{align}
Equating these for $i=j\ge2$ (by time-reversibility) and applying \Eqn{e:rj}
for $0< k<i$ (the cases $k=0$ and $k=i$ need to be handled separately),
\begin{equation} \label{e:thi}
\rr0^{i-2} \rr1^2 \bet{
  \th2 \sum_{0<k<i}\th k\th{i-k} k(i-k)-\th i\frac{i(i-1)}2 \th1^2}=0.
\end{equation}

\section{The Solutions}\label{s:2sol}
\Eqn{e:thi} holds for all $i\ge2$ if $\rr0=0$ or $\rr1=0$, leaving $\rr j=0$
by \Eqn{e:rj} for all $j\ge2$, hence $\rr0+\rr1=1$ and $\set{\th j}$ is
restricted only by the conditions $\th1>0$ and $\th+<\infty$.

\subsection{The Constant Case}\label{ss:con}
The case $\rr0=0$ leads to $\rr1=1$ and $\rr j=0$ for all $j\ne1$, so
$p(z)\equiv z$.  By \Eqn{e:phi} the joint pgf is
\[ \varphi(\zz1,\zz2,\zz3) 
     =\exp\set{P( \zz1\,\zz2\,\zz3)-\th+}, \]
so $X_1=X_2=X_3$ and all $\{X_t\}$ are identical, with an arbitrary ID
distribution.

\subsection{The IID Case}\label{ss:iid}
The case $\rr1=0$ leads to $\rr0=1$ and $\rr j=0$ for all $j\ne0$ so
$p(z)\equiv1$ and
\[ \varphi(\zz1,\zz2,\zz3) =\exp\set{P(\zz1)+P(\zz2)+P(\zz3)-3\th+} \]
by \Eqn{e:phi}, making all $\{X_t\}$ independent, with identical but
arbitrary ID distributions.

\subsection{The Poisson Case}\label{ss:poi}

Aside from these two degenerate cases, we may assume $\rr0>0$ and $\rr1>0$,
and (by \Eqn{e:rj}) rewrite \Eqn{e:thi} in the form:
\begin{align}
  \rr i&= \frac{\rr2}{\rr1^2(i-1)}\sum_{k=1}^{i-1}
          \rr k \rr{i-k},\quad i\ge2, \notag 
\intertext{whose unique solution for all integers $i\ge1$ (by induction) is} 
   \rr i&= \rr1(\rr2/\rr1)^{i-1}.\label{e:geo}
\end{align}
If $\rr2=0$, then again $\rr i=0$ for all $i\ge2$ but, by \Eqn{e:rj}, $\th
j=0$ for all $j\ge2$; thus $P(z)=\th1 z$ so each $X_t\sim\Po(\th1)$ has a
Poisson marginal distribution with mean $\th1=\th+$.  In this case
$\rr0+\rr1=1$, $p(z)=\rr0+\rr1z$, and the two-dimensional marginals (by
\Eqn{e:x1x2}) of $X_1$, $X_2$ have joint pgf
\begin{align}
    \varphi(\zz1,\zz2) 
           &= \exp\set{P\big( z\, p(\zz1)\big)-\th+
              +P(\zz1\,\rr0)-P(\rr0)}\label{e:phi-sz-p}\\
           &= \exp\set{\th1\rr0(s+z-2)+\th1\rr1(sz-1)},\notag
\end{align}
the bivariate Poisson distribution \citep [\S\thinspace37.2]
{John:Kotz:Bala:1997}, so $X_t$ is the familiar ``Poisson $\mathrm{AR}(1)$ 
Process'' of
McKenzie \citep{McKe:1985,McKe:1988} (with autocorrelation $\rho=\rr1$)
considered in \Sec {sss:pth}.  Its connection with Markov branching
processes was recognized earlier 
\citep{Steu:Verv:Wolf:1983}.
By \Eqn{e:phi-sz-p} the conditional distribution of $X_{t+1}$, given $\cF_t
:=\sigma\set{ X_s:~s\le t}$, is that of the sum of $X_t$ independent
Bernoulli random variables with pgf $p(s)$ and a Poisson innovation term with
pgf $\exp\{P(\rr0 s)-P(\rr0)\}$,
so the Markov process $X_t$ may be written recursively starting at any $t_0$ as
\begin{align*}
  X_{t_0} &\sim \Po(\th+)\\
  X_{t} &= \xi_t+\zeta_t,\text{\quad where }
    \xi_t \sim \Bi(X_{t-1}, \rr1) \text{ and }
  \zeta_t \sim \Po(\th t\rr0)\\
\end{align*}
(all independent) for $t> t_0$, the thinning construction of \Sec{s:intro}

\subsection{The Negative Binomial case}\label{ss:nb}

Finally if $\rr0>0$, $\rr1>0$, and $\rr2>0$, then (by \Eqn{e:geo}) $\rr i=
\rr1(q\rr0)^{i-1}$ for $i\ge1$ and hence (by \Eqn{e:rj}) $\th j=\alpha q^j/j$
for $j\ge1$ with $q:=(1-\rr0-\rr1)/\rr0(1-\rr0)$ and $\alpha:= \th1/q$.  The
condition $\th+<\infty$ entails $q<1$ and $\th+=-\alpha\log(1{-}q)$.  The
1-marginal distribution is $X_t\sim\NB(\alpha,p)$ with $p:=(1{-}q)$, and the
functions $P(\cdot)$ and $p(\cdot)$ are $P(z)=-\alpha \log(1-qz)$,
$p(s)=\rr0+\rr1 s/(1-q\rr0 s)$, so the joint pgf for the 2-marginal
distribution of $X_1,X_2$ is
\begin{align}
  \varphi(\zz1,\zz2)
  &= \exp\set{P\big( z\, p(\zz1)\big)-\th+ +P(\zz1\,\rr0)-P(\rr0)}\notag\\
  &=  p^{2\alpha}[(1-q\rho)-q\omrh(\zz1+\zz2)
                           +q(q-\rho)\zz1\zz2]^{-\alpha}\label{e:phi-sz-nbb}
\end{align}
with one-step autocorrelation $\rho:=(1{-}\rr0)^2/\rr1$.  This bivariate
distribution was introduced as the ``compound correlated bivariate
Poisson''\citep{Edwa:Gurl:1961}, but we prefer to call it the Branching
Negative Binomial distribution.  In the branching formulation $X_{t}$ may
be viewed as the sum of $X_{t-1}$ iid random variables with pgf $p(s)=\rr0
+\rr1 s/(1-q\rr0 s)$ and one with pgf $\exp\set{P(s\rr0)-P(\rr0)}
=(1-q\rr0)^\alpha (1-q\rr0\,s)^{-\alpha}$.  The first of these may be
viewed as $Y_t$ plus a random variable with the $\NB(Y_t,1{-}q\rr0)$
distribution, for $Y_t\sim\Bi(X_{t-1},1-\rr0)$, and the second has the
$\NB(\alpha,1{-}q\rr0)$ distribution, so a recursive updating scheme
beginning with $X_{t_0} \sim \NB(\alpha,p)$ is:
\[ X_t = Y_t+\zeta_t,\text{\quad where }
     Y_t \sim \Bi(X_{t-1}, ~1{-}\rr0)\text{ and }
         \zeta_t\sim \NB(\alpha+Y_t,~1{-}q\rr0).
\]
In the special case of $\rho=q$ the joint pgf simplifies to $\varphi
(\zz1,\zz2) = p^\alpha[1+q(1-\zz1-\zz2)]^{-\alpha}$ and the joint
distribution of $X_1,X_2$ reduces to the negative trinomial distribution
\citep [Ch.~36] {John:Kotz:Bala:1997} with pmf
\begin{align*}
\P[X_1=i,X_2=j] &= 
  \frac {\Gamma(\alpha+i+j)} {\Gamma(\alpha)~i!~j!} 
  \cet{\frac{1-q}{1+q}}^\alpha\,\cet{\frac{q}{1+q}}^{i+j}
\end{align*}
and simple recursion $X_t\mid X_{t-1}\sim\NB\big(\alpha+X_{t-1},
~\frac1{1+q}\big)$.

\subsection{Results}\label{ss:thm}

We have just proved:
\begin{thm}\label{t:struc}
  Let $\set{X_t}$ be a Markov process indexed by $t\in\bbZ$ taking values in
  the non-negative integers $\bbN_0$ that is stationary, time-reversible, has
  infinitely-divisible marginal distributions of all finite orders, and
  satisfies $\P[X_t=1]>0$.  Then $\set{X_t}$ is one of four processes:
\begin{enumerate}
\item \label{i:const}
  $X_t\equiv X_0\sim\mu_0(dx)$ for an arbitrary ID distribution $\mu_0$
  on $\bbN_0$ with $\mu_0(\{1\})>0$;
\item \label{i:iid}
   $X_t\iid \mu_0(dx)$ for an arbitrary ID distribution $\mu_0$ on
  $\bbN_0$ with $\mu_0(\{1\})>0$;
\item \label{i:po}
  For some $\theta>0$ and $0<\rho<1$,
  $X_t\sim\Po(\theta)$ with bivariate joint generating function
  \[
   \E\left[\zz1^{X_1}~\zz2^{X_2}\right]
        = \exp\set{\theta \omrh(\zz1-1) +
                   \theta \omrh(\zz2-1) +
                   \theta  \rho(\zz1\zz2-1) }
  \]
  and hence correlation $\Corr(X_s,X_t)= \rho^{|s-t|}$ and recursive update
\[ X_t = \xi_t+\zeta_t,\text{\quad where }
     \xi_t \sim \Bi(X_{t-1}, ~\rho)\text{ and }
     \zeta_t\sim \Po\big(\theta\omrh);
\]
\item \label{i:nb}
  For some $\alpha>0$, $0<p<1$, and $0<\rho<1$, $X_t\sim\NB(\alpha,p)$, with
  bivariate joint generating function
\[
   \E\left[\zz1^{X_1}~\zz2^{X_2}\right] =
    p^{2\alpha}[(1-q\rho)-q\omrh(s+z)+q(q-\rho)sz]^{-\alpha}
\]
where $q=1{-}p$, and hence correlation $\Corr(X_s,X_t)= \rho^{|s-t|}$ and
recursive update
\[ X_t = Y_t+\zeta_t,\text{\quad where }
     Y_t \sim \Bi\big(X_{t-1}, ~\rho\,p/(1-\rho q)\big)\text{ and }
     \zeta_t\sim \NB\big(\alpha+Y_t,~p/(1-\rho q)\big).
\]
\end{enumerate}
\end{thm}
Note the limiting cases of autocorrelation $\rho=1$ and $\rho=0$ in cases
\ref{i:po}., \ref{i:nb}. are subsumed by the degenerate cases
\ref{i:const}. and \ref{i:iid}., respectively.  The theorem follows from
this.

From this theorem follows:

\begin{thm}\label{t:thin}
  Let $\set{\mu^\theta:~\theta\ge0}$ be an ID semigroup of probability
  distributions on the nonnegative integers $\bbN_0$ with $\mu^\theta
  (\{1\})>0$.  Fix $\theta>0$ and $0<\rho<1$ and let $\set{X_t}$ be the
  ``thinning process'' of \Eqn{e:thin} in \Sec{ss:thin} with the
  representation
\begin{equation}\label{e:rev}
    X_{t-1}=\xi_t+\eta_t \qquad X_t=\xi_t+\zeta_t 
\end{equation}
  for each $t\in\bbZ$ with independent
   \[ \xi_t\sim\mu^{\rho\theta}(d\xi)\qquad
      \eta_t\sim\mu^{\omrh\theta}(d\eta)\qquad
      \zeta_t\sim\mu^{\omrh\theta}(d\zeta).\]
   Then $X_t$ is Markov, stationary, time-reversible, and nonnegative integer
   valued, but it does not have infinitely-divisible marginal distributions
   of all orders unless $\{\mu^\theta\}$ is the Poisson family.
\end{thm}
\proof By construction $X_t$ is obviously Markov and stationary.  The joint
distribution of the process at consecutive times is symmetric (see \Eqn {e:rev})
since the marginal and conditional pmfs
\[ p(x) := \mu^\theta(\{x\}),\qquad  q(y\mid x) := 
       \frac{\sum_z \mu^{\rho\theta}(\{z\})
                   ~\mu^{\omrh\theta}(\{x-z\})
                   ~\mu^{\omrh\theta}(\{y-z\})}
            {\mu^{\theta}(\{x\})}
\]
of $X_t$ and $X_t\mid X_{t-1}$ satisfy the symmetric relation
\[ p(x)~q(y\mid x) = q(x\mid y)~ p(y). \] 
Applying this inductively, for any $s<t$ and any $\{x_s,\cdots,x_t\}\subset
\bbN_0$ we find
\begin{align*}
\P[X_s=x_s,\cdots,X_t=x_t] 
   &= p(x_s)\hspace{4mm}  q(x_{s+1}\mid x_s) \hspace{1mm} 
          q(x_{s+2}\mid x_{s+1})\cdots q(x_t\mid x_{t-1})\\
   &= q(x_s\mid x_{s+1}) p(x_{s+1}) \hspace{1mm} 
          q(x_{s+2}\mid x_{s+1}) \cdots q(x_t\mid  x_{t-1})\\
   &=\cdots\\
   &= q(x_s\mid x_{s+1}) q(x_{s+1}\mid x_{s+2})\cdots
                 q(x_{t-1}\mid x_t)\hspace{1mm} p(x_t),
\end{align*}
and so the distribution of $X_t$ is time-reversible.  Now suppose that it is
also ID.  Then by \Thm{t:struc} it must be one of the four specified
processes: constant, iid, branching Poisson, or branching negative binomial.

Since $\rho<1$ it cannot be the constant $\set{X_t\equiv X_0}$ process; since
$\rho>0$ it cannot be the independent $\set{X_t\iid\mu^\theta(dx)}$ process.
The joint generating function $\phi(s,z)$ at two consecutive times for the
negative binomial thinning process, given in \Eqn {e:phi-sz-nbt}, differs
from that for the negative binomial branching process, given in \Eqn
{e:phi-sz-nbb}.  The only remaining option is the Poisson branching process
of \Sec {sss:pth}.
\endproof

\begin{thm}\label{t:meas}
  Let $\set{\mu^\theta:~\theta\ge0}$ be an ID semigroup of probability
  distributions on the nonnegative integers $\bbN_0$ with $\mu^\theta
  (\{1\})>0$.  Fix $\theta>0$ and $0<\rho<1$ and let $\set{X_t}$ be the
  ``random measure process'' of \Eqn{e:meas} in \Sec{ss:meas}.  Then $X_t$ is
  ID, stationary, time-reversible, and nonnegative integer valued, but it is
  not a Markov process unless $\{\mu^\theta\}$ is the Poisson family.
\end{thm}
\proof By construction $X_t$ is ID, stationary, and time-reversible; suppose
that it is also Markov.  Then by \Thm{t:struc} it must be one of the four
specified processes: constant, iid, branching Poisson, or branching negative
binomial.

Since $\rho<1$ it cannot be the constant $\set{X_t\equiv X_0}$ process; since
$\rho>0$ it cannot be the independent $\set{X_t\iid\mu^\theta(dx)}$ process.
The joint generating function $\phi(s,z)$ at two consecutive times for the
negative binomial random measure process coincides with that for the negative
binomial thinning process, given in \Eqn {e:phi-sz-nbt}, and differs from
that for the negative binomial branching process, given in \Eqn
{e:phi-sz-nbb}.  The only remaining option is the Poisson branching process
of \Sec {sss:pth}.
\endproof

\section{Continuous Time}\label{s:cont}
Now consider $\bbN_0$-valued time-reversible stationary Markov processes
indexed by continuous time $t\in\bbR$.  The restriction of any such process
to $t\in\bbZ$ will still be Markov, hence MISTI, so there can be at most two
non-trivial ones--- one with univariate Poisson marginal distributions, and
one with univariate Negative Binomial distributions.  Both do in fact exist.

\subsection{Continuous-Time Poisson Branching Process}\label{ss:ctpoi}

Fix $\theta>0$ and $\lambda>0$ and construct a nonnegative integer-valued
Markov process with generator
\begin{subequations}\label{e:GenPoi}
\begin{align}
\fA f(x) &= \frac\partial{\partial s}
                 \E[ f(X_s)-f(X_t)\mid X_t=x ] \Big|_{s=t}\notag\\
         &= \lambda\theta \big[f(x+1)-f(x)]\big]
          + \lambda x \big[f(x-1)-f(x)\big]\label{e:GenPoi.a}
\intertext{or, less precisely but more intuitively, for all $i,j\in\bbN_0$
  and $\eps>0$,}
\P\big[X_{t+\eps}=i\mid X_t=j\big] &= o(\eps) +
   \begin{cases}
    \eps\lambda\theta       &i=j+1\\
    1-\eps\lambda(\theta+j) &i=j\\
    \eps\lambda j           &i=j-1\\
    \end{cases}\label{e:GenPoi.b}
\end{align}
\end{subequations}
$X_t$ could be described as a linear death process with immigration.  In
\Sec{ss:marg} we verify that its univariate marginal distribution and
autocorrelation are
\begin{align*}
            X_t& \sim\Po(\theta)\\
 \Corr(X_s,X_t) &= e^{-\lambda|s-t|},
\end{align*}
and its restriction to integer times $t\in\bbZ$ is precisely the process
described in \Sec{s:2sol} item 3, with one-step autocorrelation
$\rho=e^{-\lambda}$.

\subsection{Continuous-Time Negative Binomial Branching Process}\label{ss:ctnb}

Now fix $\theta>0$, $\lambda>0$, and $0<p<1$ and construct a nonnegative
integer-valued Markov process with generator
\begin{subequations}\label{e:GenNB}
\begin{align}
\fA f(x) &= \frac\partial{\partial s}
                 \E[ f(X_s)-f(X_t)\mid X_t=x ] \Big|_{s=t}\notag\\
         &= \frac{\lambda(\alpha+x)\omp}p \big[f(x+1)-f(x)]\big]
          + \frac{\lambda x}p \big[f(x-1)-f(x)\big]\label{e:GenNB.a}
\intertext{or, for all $i,j\in\bbN_0$ and $\eps>0$,}
\P\big[X_{t+\eps}=i\mid X_t=j\big] &= o(\eps) +
   \begin{cases}
    \eps\lambda(\alpha+j)\omp/p       &i=j+1\\
    1-\eps\lambda[(\alpha+j)\omp+j]/p &i=j\\
    \eps \lambda j/p        &i=j-1,
    \end{cases}\label{e:GenNB.b}
\end{align}
\end{subequations}
so $X_t$ is a linear birth-death process with immigration.  The univariate
marginal distribution and autocorrelation (see \Sec{ss:marg}) are now
\begin{align*}
            X_t& \sim \NB(\alpha,p)\\
 \Corr(X_s,X_t) &= e^{-\lambda|s-t|},
\end{align*}
and its restriction to integer times $t\in\bbZ$ is precisely the process
described in \Sec{s:2sol} item 4, with autocorrelation $\rho =e^{-\lambda}$.

\subsection{Markov Branching (Linear Birth/Death) Processes }\label{ss:pops}

The process $X_t$ of \Sec{ss:ctpoi} can also be described as the size of a
population at time $t$ if individuals arrive in a Poisson stream with rate
$\lambda\theta$ and die or depart independently after exponential holding
times with rate $\lambda$; as such, it is a continuous-time Markov branching
process.

Similarly, that of \Sec{ss:ctnb} can be described as the size of a population
at time $t$ if individuals arrive in a Poisson stream with rate
$\lambda\alpha \omp/p$, give birth (introducing one new individual)
independently at rate $\lambda \omp/p$, and die or depart at rate
$\lambda/p$.  In the limit as $p\to1$ and $\alpha\to\infty$ with
$\alpha\omp\to\theta$ this will converge in distribution to the Poisson
example of \Sec{ss:ctpoi}.

\subsection{Marginal Distributions}\label{ss:marg}
Here we verify that the Poisson and Negative Binomial distributions are the
univariate marginal stationary distributions for the Markov chains with
generators $\fA$ given in
\Eqn{e:GenPoi} and \Eqn{e:GenNB}, respectively.

Denote by $\pi^0_i=\P[X_t=i]$ the pmf for $X_t$ and by
$\pi^\eps_i=\P[X_{t+\eps}=i]$ that for $X_{t+\eps}$, and by
$\varphi_0(s)=\E[s^{X_t}]$ and $\varphi_\eps(s)=\E[s^{X_{t+\eps}}]$ their
generating functions.  The stationarity requirement that $\varphi_0(s)\equiv
\varphi_\eps(s)$ will determine $\varphi(s)$ and hence $\{\pi_i\}$ uniquely.
\subsubsection{Poisson}\label{sss:ctPo}
From \Eqn{e:GenPoi.b} for $\eps>0$ we have
\begin{alignat}5
\pi^\eps_i  &= \eps\lambda\theta \pi^0_{i-1}&
         &+ [1-\eps\lambda(\theta+i)] \pi^0_{i}&
         &+\eps\lambda(i+1) \pi^0_{i+1}&
         &+ o(\eps). \notag
\intertext{Multiplying by $s^i$ and summing, we get:}
\varphi_\eps(s)
         &= \eps\lambda\theta s \sum_{i\ge1} s^{i-1} \pi^0_{i-1}&
         &+ [1-\eps\lambda\theta] \varphi_0(s)
          - \eps\lambda s \sum\iz i s^{i-1} \pi^0_{i}&
         &+ \eps\lambda \sum\iz (i+1)s^i \pi^0_{i+1}&
         &+ o(\eps) \notag\\
         &= \eps\lambda\theta s \varphi_0(s) &
         &+ [1-\eps\lambda\theta] \varphi_0(s)
          - \eps\lambda s \varphi_0'(s)&
         &+ \eps\lambda \varphi_0'(s) &
         &+ o(\eps) \notag
\end{alignat}
so
\begin{align*}
\varphi_\eps(s)-\varphi_0(s)
         &= \eps\lambda(s-1)\bet{\theta\varphi_0(s)-\varphi_0'(s)}+ o(\eps)
         \notag 
\intertext{and stationarity ($\varphi_0(s)\equiv\varphi_\eps(s)$) entails
  $\lambda=0$ or $\varphi_0'(s)/\varphi_0(s) \equiv \theta$, so
  $\log\varphi_0(s) \equiv (s-1)\theta$ and:}
          \varphi_0(s) &= \exp\set{(s-1)\theta}\notag\\
\end{align*}
so $X_t\sim\Po(\theta)$ is the unique stationary distribution.

\subsubsection{Negative Binomial}\label{sss:ctNB}
From \Eqn{e:GenNB.b} for $\eps>0$ we have
\begin{align*}
\pi^\eps_i  &= (\eps\lambda\omp/p)(\alpha+i-1)~\pi^0_{i-1}
         + \{1-(\eps\lambda/p)[(\alpha+i)\omp+i]\}~ \pi^0_{i}
         + (\eps\lambda/p)(i+1)~ \pi^0_{i+1}+ o(\eps) \notag\\
\varphi_\eps(s)
         &= (\eps\lambda\omp/p)\alpha~ s\varphi_0(s) + 
            (\eps\lambda\omp/p)~ s^2\varphi_0'(s)\\
         &+ \varphi_0(s) - (\eps\lambda\omp/p) \alpha\varphi_0(s)
          - (\eps\lambda/p)(\omp+1)~s\varphi_0'(s)\\
         &+ (\eps\lambda/p)~\varphi_0'(s)+ o(\eps) \notag\\
\varphi_\eps(s)-\varphi_0(s)
         &=(\eps\lambda/p)~\set{
            \varphi_0(s)~\alpha\omp (s-1)
          + \varphi_0'(s)~[\omp s^2 -(\omp+1)s +1]}+ o(\eps) \notag\\
         &=(\eps\lambda/p)(s-1)~\set{
            \varphi_0(s)~\alpha\omp
          + \varphi_0'(s)~(\omp s-1)}+ o(\eps) \notag
\end{align*}
so either $\lambda=0$ (the trivial case where $X_t\equiv X$) or $\lambda>0$ and:
\begin{align*}
\varphi_0'(s)/\varphi_0(s) &= \alpha\omp (1-\omp s)^{-1}\\
         \log \varphi_0(s) &= -\alpha\log (1-\omp s)+\alpha\log(p)\\
              \varphi_0(s) &= p^\alpha (1-\omp s)^{-\alpha}
\end{align*}
and $X_t\sim\NB(\alpha,p)$ is the unique stationary distribution.

\subsubsection{Alternate Proof}\label{sss:alt}
A detailed-balance argument \citep[\pg105]{Hoel:Port:Ston:1972} shows that
the stationary distribution $\pi_i:=\P[X_t=i]$ for linear birth/death chains
is proportional to
\begin{align*}
  \pi_i &\propto \prod_{0\le j<i} \frac{\beta_{j}}{\delta_{j+1}}\\
  \intertext{where $\beta_j$ and $\delta_j$ are the birth and death rates
    when $X_t=j$, respectively.  For the Poisson case, from \Eqn{e:GenPoi.b}
    this is} \pi_i &\propto \prod_{0\le j<i} \frac{\lambda\theta}{\lambda
    (j+1)} = \theta^i/i!, \intertext{so $X_t\sim\Po(\theta)$, while for the
    Negative Binomial case from \Eqn{e:GenNB.b} we have} \pi_i &\propto
  \prod_{0\le j<i} \frac{\lambda(\alpha+j)\omp/p} {\lambda (j+1)/p} =
  \frac{\Gamma(\alpha+i)}{\Gamma(\alpha)\,i!}~ \omp^i,
\end{align*}
so $X_t\sim\NB(\alpha,p)$.  In each case the proportionality constant is
$\pi_0=P[X_t=0]$:  $\pi_0=e^{-\theta}$ for the Poisson case, and
$\pi_0=p^\alpha$ for the negative binomial.

\subsubsection{Autocorrelation}\label{sss:auto}
Aside from the two trivial (iid and constant) cases, MISTI processes have
finite $p$th moments for all $p<\infty$ and, in particular, have finite
variance and well-defined autocorrelation.  It follows by the Markov property
and induction that the autocorrelation must be of the form
\[ \Corr[X_s,X_t] = \rho^{-|t-s|} \] for some $\rho\in[-1,1]$.  In both the
Poisson and negative binomial cases the one-step autocorrelation $\rho$ is
nonnegative; without loss of generality we may take $0<\rho<1$.

\hide{
    From \Eqn{e:GenNB.b} the conditional expectation of $X_{t+\eps}$ for $\eps>0$
    given $\cF_t$ is
    \begin{align*}
    \E[X_{t+\eps} \mid X_t] &= X_t(1-\eps\lambda) 
                                    + \eps\lambda\alpha\omp/p +o(\eps)
    \intertext{and hence for $t>0$ and $\eps>0$ the autocorrelation
    $\rho(t)=\Corr(X_0,X_t)$ satisfies} 
    \rho(t+\eps)&= \Corr[X_0,X_{t+\eps}]\\
            &=  \rho(t) [1-\eps\lambda] + o(\eps)
    \intertext{so $\rho'(t) = -\lambda\rho(t)$, $\rho(t)=e^{-\lambda t}$ for
      $t\ge0$ and, by symmetry,}
       \Corr[X_s, X_t] &=  e^{-\lambda|t-s|}
    \end{align*}
    as claimed.  A similar argument beginning with \Eqn{e:GenPoi.b} shows
    $\E[X_{t+\eps}\mid X_t] = X_t(1-\eps\lambda) + \eps\lambda\theta+o(\eps)$ and
    hence $\rho(t)=e^{-\lambda|t|}$ for the Poisson example of
    \Sec{sss:ctPo}.
}

\section{Discussion}\label{s:disc}

The condition $\mu^\theta(\{1\})>0$ introduced in \Sec{ss:3dim} to avoid
trivial technicalities is equivalent to a requirement that the support
$\spt(\mu^\theta) =\bbN_0$ be all of the nonnegative integers.  Without this
condition, for any MISTI process $X_t$ and any integer $k\in\bbN$ the process
$Y_t=k\,X_t$ would also be MISTI, leading to a wide range of essentially
equivalent processes.

The branching approach of \Sec {ss:pops} could be used to generate a wider
class of continuous-time stationary Markov processes with ID marginal
distributions \citep {Verv:1979, Steu:Verv:Wolf:1983}.  If families of size
$k\ge1$ immigrate independently in Poisson streams at rate $\lambda_k$, with
$\sum_{k\ge1} \lambda_k\log k<\infty$, and if individuals (after independent
exponential waiting times) either die (at rate $\delta>0$) or give birth to
some number $j\ge1$ of progeny (at rate $\beta_j\ge0$), respectively, with
$\delta> \sum_{j\ge1} j\,\beta_j$, then the population size $X_t$ at time $t$
will be a Markov, infinitely-divisible, stationary processes with nonnegative
integer values.  Unlike the MISTI processes, these may have infinite $p$th
moments if $\sum_{k\ge1} \lambda_k k^p=\infty$ for some $p>0$ and, in
particular, may not have finite means, variances, or autocorrelations.

Unless $\lambda_k =0$ and $\beta_j=0$ for all $k,j>1$, however, these will
not be time-reversible, and hence not MISTI.  Decreases in population size
are always of unit size (necessary for the Markov property to hold), while
increases might be of size $k>1$ (if immigrating family sizes exceed one) or
$j>1$ (if multiple births occur).

\section*{Acknowledgments}
The authors would like to thank
Xuefeng Li,
Avi     Mandelbaum,
Yosef   Rinott,
Larry   Shepp,
and 
Henry   Wynn
for helpful conversations.
 This work was supported in part by National
 Science Foundation grants 
 DMS--1228317 and 
 DMS-2015382      
 and National Air and Space Administration Applied Information Science
 Research Program grant NNX09AK60G.  
 Larry Brown is sorely missed, both for his deep intellect and his
 delightful charm, generosity, and humanity.

\hide{
  \section{Leftovers}\label{s:left}
  \begin{itemize}
  \item The requirement $\th1>0$ (made just above \Eqn{e:pP}) probably isn't
  necessary... if not, set $m\equiv\min\{i:~\th i>0\}$ and define $\rr j :=
  \lll jm+/\th m$ and I expect we can still characterize all the MISTI's.  My
  guess is that they will be the same two degenerate ones, plus $m$ times the
  Po$^*$ and $\mathsf{NT}(\alpha,p)$ solutions (for example, if $\th1=0$ and
  $\th2>0$, then I expect all $X_t$ to be even, and for $X_t/2$ to be one of
  our solutions), but I haven't checked that.  The key will be showing that
  if $\th m>0$ and $\th n>0$ then $\E[\zz1^{X_1}\,\zz3^{X_3}\mid X_2=
  \mathsf{lcm}(m,n)]$ factors properly in two ways---  so that the
  conditional distribution of $X_3$ given $X_2=j$ won't depend on what mix
  of $m$'s and $n$'s comprise $j$, and hence won't depend on $X_1$.
  \end{itemize}
}
%
%
\ifBib
\bibliography{arp}
\else
\newcommand{\noopsort}[1]{}

\fi 
\par\bigskip
\centerline{\begin{tabular}{l@{\qquad}l}
Robert L. Wolpert            &Lawrence D. Brown\\
Department of Statistical Science &Department of Statistics\\
Duke University              &Wharton School, University of Pennsylvania\\
Durham, NC 27708-0251 USA    &Philadelphia, PA 19104 USA\\
              \texttt{rlw@duke.edu}&(Deceased)\\ 
\url{http://www.stat.duke.edu/~rlw/}\\ 
\end{tabular}}

\end{document}